\documentclass[12pt]{article}

\usepackage{amssymb, amsmath}

\newcommand{\eop}{\bigstar}  

\newcommand{\card}[1]{{\vert #1 \vert} }

\newenvironment{Proof}{\noindent{\bf Proof.}}{\par\bigskip} 

\newenvironment{sadrzaj}{\centerline{\bf Sadr\v zaj}}{\par\bigskip}

\newtheorem{THEOREM}{Theorem}[section]

\newtheorem{Conclusion}[THEOREM]{Conclusion}

\newtheorem{LEMMA}[THEOREM]{Lemma}

\newtheorem{Main Theorem}[THEOREM]{Main Theorem}
\newenvironment{main Theorem}{\begin{Main Theorem}} 
{\end{Main Theorem}}

\newtheorem{Theorem}[THEOREM]{Theorem}

\newtheorem{Definition}[THEOREM]{Definition}

\newtheorem{Conventions}[THEOREM]{Conventions}

\newtheorem{Main Definition}[THEOREM]{Main Definition}
\newenvironment{main definition}{\begin{Main Definition}}
{\end{Main Definition}}

\newtheorem{Lemma}[THEOREM]{Lemma}

\newtheorem{Notation}[THEOREM]{Notation}

\newtheorem{Note}[THEOREM]{Note}

\newtheorem{Observation}[THEOREM]{Observation}

\newtheorem{Remark}[THEOREM]{Remark}

\newtheorem{Main Fact}[THEOREM]{Main Fact}
\newenvironment{main Fact}{\begin{Main Fact}}{\end{Main Fact}}

\newtheorem{Fact}[THEOREM]{Fact}

\newtheorem{Subfact}[THEOREM]{Subfact}

\newtheorem{Claim}[THEOREM]{Claim}

\newtheorem{Main Claim}[THEOREM]{Main Claim}
\newenvironment{main claim}{\begin{Main Claim}}{\end{Main Claim}}

\newtheorem{Corrolary}[THEOREM]{Corrolary}

\newtheorem{Subclaim}[THEOREM]{Subclaim}

\newtheorem{Corollary}[THEOREM]{Corollary}

\newtheorem{Proposition}[THEOREM]{Proposition}

\newtheorem{Discussion}[THEOREM]{Discussion}

\newenvironment{Proof of the Subfact}
{\noindent{\bf Proof of the Subfact.}}{\par\bigskip}

\newenvironment{Proof of the Fact}
{\noindent{\bf Proof of the Fact.}}{\par\bigskip}

\newenvironment{Proof of the Lemma}
{\noindent{\bf Proof of the Lemma.}}{\par\bigskip}

\newenvironment{Proof of the Claim}
{\noindent{\bf Proof of the Claim.}}{\par\bigskip}

\newenvironment{Proof of the Subclaim}
{\noindent{\bf Proof of the Subclaim.}}{\par\medskip}

\newenvironment{Proof of the Main Claim}
{\noindent{\bf Proof of the Main Claim.}}{\par\bigskip}

\catcode`\@=11
\def\@begintheorem#1#2{\rm \trivlist \item[\hskip \labelsep{\bf #1\
#2.}]}
\def\@opargbegintheorem#1#2#3{\rm \trivlist
      \item[\hskip \labelsep{\bf #1\ #2\ (#3).}]}
\catcode`\@=12

\newcommand{\into}{\rightarrow}

\newcommand{\deq}{\buildrel{\rm def}\over =}

\newcommand{\DD}{{\cal D}}

\newcommand{\PP}{{\cal P}}

\newcommand{\UU}{{\cal U}}

\title{A Note on the Splitting Property in Strongly Dense Posets of
Size $\aleph_0$}

\author{Mirna D\v zamonja\\
School of Mathematics\\
University of East Anglia\\
Norwich, NR4 7TJ, UK\\
{\scriptsize{M.Dzamonja@uea.ac.uk}}}

\date{published in Radovi Matemati{\v c}ki, 8, no. 2 (1998), pg. 321-326}

\begin{document}
\maketitle

\begin{abstract}
We show that it is not true that every countable
infinite strongly dense poset has the splitting property, so
answering a question of R. Ahlswede, P.L. Erd\"os and N. Graham.
\footnote{We thank P\'eter L. Erd\"os and Niall Graham for
discussing the problem, and Mary-Ellen Rudin and Jean-Marc 
Vanden-Broeck for listening to the proof.}
\end{abstract}

\vfill

\baselineskip=16pt
\binoppenalty=10000
\relpenalty=10000
\raggedbottom

\section{Introduction} 
\begin{Definition} Suppose that $\PP=\langle P,\le_P=\le\rangle$
is a partially ordered set, and $H\subseteq P$. We define

(1)
\[
\DD(H)\deq\{x\in P:\,(\exists s\in H)\,(x\le s)\}
\]
and 
\[
\UU(H)\deq\{x\in P:\,(\exists s\in H)\,(x\ge s)\}.
\]

(2) $H$ is called an {\em antichain} or a {\em Sperner system} of $P$ if no
two elements of $H$ are comparable.

(3) An antichain $A$ of $P$ is {\em maximal} iff 
\[
[P\supseteq A'\supseteq A\,\,\&\,\,A'\mbox{ an antichain of }P]\implies
A'=A.
\]

(4) An antichain $A$ of $P$ is said to {\em split} iff there is
a partition $\{D,U\}$ of $A$ such that $P=\DD(D)\cup\UU(U)$.

(5) $P$ is said to have the {\em splitting property} iff every maximal antichain
of $P$ splits.
\end{Definition}

\begin{Notation} In a poset $\PP=\langle P,\le_P\rangle$:

(1) For $x,y\in P$ we let $x<_P y$ iff $x\le_P y$ and $x$ and $y$ are
not equal.

(2) For $x,y\in P$
\[
(x,y)\deq\{z\in P:\,x<_Pz<_P y\},
\]
and sets of the form $(x,y)$ for $x<_Py\in P$ are called open intervals
of $P$.

(3) If $x,y\in P$ are incomparable, i.e. neither $x\le_P y$
nor $y\le_P x$, we denote this by $x><y$.
\end{Notation}

\begin{Definition} A poset $\PP=\langle P,\le_P\rangle$
is said to be {\em strongly dense} iff for every non-empty
open interval $(x,y)$ of $P$, there are $z_1,z_2\in (x,y)$
such that $z_1><z_2$.
\end{Definition}

In \cite{AEG}, it was shown that every finite strongly dense
poset has the splitting property. It was conjectured that
the same is true for countably infinite strongly dense
posets. We give an example
showing that the conjecture is not true. In \S\ref{theorem}
we give a formal construction of the example and the proof that
it is as required. The idea is however very simple: take the
binary tree ${}^{<\omega}2$. Consider the maximal antichain
$A=\{\langle 0\rangle,\langle 1\rangle\}$ in the tree: it does not split.
But the binary tree is not strongly dense. To remedy this, we replace
every node of the tree by a copy of ${}^{<\omega}2$, and continue
$\omega$ many times, extending the partial order
in the natural way as we go. Taking the union at the end, we arrive
at a strongly dense partial order, and the version of $A$
on the second level of the construction can be
used to exemplify a maximal antichain which does not split.

The splitting property for arbitrary posets is further discussed
in \cite{E} where a sufficient condition for the splitting
is given. One can consult \cite{E}
or \cite{AEG} for more on the history
of the problem and further references.

\section{The example}\label{theorem}

\begin{Notation} (1) ${}^{<\omega}2\deq\bigcup_{n<\omega}
\{s|\,s:n\into\{0,1\}\mbox{ is a function}\}$, ordered 
by being an initial segment.

(2) For $s,t\in {}^{<\omega}2$ we let $s\unlhd t$ iff $s$ is an initial
segment of $t$.

(3) $\langle\rangle$ stands for the empty sequence.

\end{Notation}

\begin{Theorem}\label{example} There is a countably infinite
strongly dense poset $\langle C,\le\rangle$
which does not have the splitting property.

(Moreover, the poset $\langle C,\le\rangle$ has the property described in
the statement of Claim \ref{one}(3) below). 
\end{Theorem}

\begin{Proof} We start by defining the poset:

\begin{Definition}\label{main}
(1) By induction on $n<\omega$ we define $C_n$:

\underline{$n=0$}.
\[
C_0\deq\{\langle\emptyset,s,0\rangle:\,s\in{}^{<\omega}2\}
\]
ordered by $\le_0$ given by
\[
\langle\emptyset,s,0\rangle\le_0\langle\emptyset,t,0\rangle
\iff s\unlhd t.
\]

\underline{$n+1$}. Given $C_n$ and a partial order $\le_n$ on $C_n$.
We define:
\[
C_{n+1}\deq C_n\bigcup\left\{\langle x,s,n+1\rangle:\,s\in {}^{<\omega}2
\setminus\{\langle\rangle\}\,\,\&\,\,x\in C_n\right\},
\]
ordered by letting $\langle a_1,b_1,c_1\rangle\le_{n+1}
\langle a_2,b_2,c_2\rangle$
iff one of the following occurs:
\begin{description}
\item{(i)} $c_1,c_2\le n$ and $\langle a_1,b_1,c_1\rangle
\le_n\langle a_2,b_2,c_2\rangle$,
\item{(ii)} $c_1\le n\,\,\&\,\,c_2=n+1\,\,\&\,\,\langle a_1,b_1,c_1\rangle
\le_n a_2$,
\item{(iii)} $c_1=n+1\,\,\&\,\, c_2\le n \,\,\&\,\,
a_1<_n \langle a_2,b_2,c_2\rangle$,
\item{(iv)} $c_1=c_2=n+1\,\,\&\,\,a_1=a_2\,\,\&\,\,b_1\unlhd b_2$,
\item{(v)} $c_1=c_2=n+1\,\,\&\,\,a_1<_n a_2$.
\end{description}

(2) We let $C\deq\bigcup_{n<\omega}C_n$, ordered by
\[
x\le y\iff x\le_n y\mbox{ for the minimal }n \mbox{ such that }
x,y\in C_n.
\]

\end{Definition}

\begin{Remark}\label{agree} Observe that in Definition \ref{main}(2),
for $x,y\in C$
we have 
\[
x\le y\iff x\le_n y\mbox{ for any } n\mbox{ such that }x,y\in C_n,
\]
because of the item (i) in Definition \ref{main}(1).
${\eop}_{\ref{agree}}$
\end{Remark}

\begin{Claim}\label{one} (1) $C$ is countable.

(2) $\le$ is a partial order on $C$.

(3) $\langle C,\le \rangle$ is strongly dense, moreover for
every $x,y\in P$
\[
x<y\implies(\exists Z\subseteq (x,y))\,[\card{Z}=\aleph_0\,\,\&\,\,
z_1\neq z_2\implies z_1><z_2].
\]
\end{Claim}

\begin{Proof of the Claim} (1) Obvious, proving by induction on $n$
that $C_n$ is countable.

(2) The only part which needs checking is that $\le$ is transitive.
As $\le=\bigcup_{n<\omega}\le_n$
is an increasing union, it suffices to show that every $\le_n$
is transitive. We do this by induction on $n$.

\underline{$n=0$}. Obvious.

\underline{$n+1$}. Given $\langle a_l, b_l,c_l\rangle$ for $l\in 
\{0,1,2\}$ such that
\[
\langle a_l, b_l,c_l\rangle\le_{n+1}
\langle a_{l+1}, b_{l+1},c_{l+1}\rangle \mbox{ for }l\in\{0,1\}.
\]
We shall prove that $\langle a_0, b_0,c_0\rangle
\le_{n+1}\langle a_2, b_2,c_2\rangle$
by considering several cases.

\smallskip

\underline{Case 1}. $c_l\le n$ for $l\in 
\{0,1,2\}$.

Follows by the induction hypothesis.

\underline{Case 2}. $c_0\le n$ and $c_1=c_2=n+1$.

Hence we have $\langle a_0,b_0,c_0\rangle\le_n a_1\le_n a_2$
(by cases (ii), (iv), (v) of Definition \ref{main}(1)),
so the conclusion follows by case (ii) of Definition \ref{main}(1)
and the induction hypothesis.

\underline{Case 3}. $c_0, c_2\le n$ and $c_1=n+1$.

Hence $\langle a_0,b_0,c_0\rangle\le_n a_1 <_n \langle a_2,b_2,c_2
\rangle$, (by cases (ii) and (iii) of Definition \ref{main}(1)),
so the conclusion follows by the induction hypothesis.

\underline{Case 4}. $c_0=n+1$ and $c_1,c_2\le n$.

Hence $a_0<_n\langle a_1,b_1,c_1\rangle \le_n \langle a_2,b_2,c_2
\rangle$, (by (ii), (iv) and (v) of Definition \ref{main}(1)),
so the conclusion follows by case (iii) of
Definition \ref{main}(1).

\underline{Case 5}. $c_0=c_2=n+1$ and $c_1\le n$.

Hence $a_0<_n\langle a_1,b_1,c_1\rangle \le_n a_2
$, (by (iii) and (ii) of Definition \ref{main}(1)),
so the conclusion follows by case (v) of
Definition \ref{main}(1).

\underline{Case 6}. $c_0,c_1\le n$ and $c_2= n+1$.

Hence $\langle a_0,b_0,c_0\rangle
\le_n\langle a_1,b_1,c_1\rangle \le_n a_3$,
(by (i) and (ii) of Definition \ref{main}(1)),
so the conclusion follows by case (ii) of
Definition \ref{main}(1).

\underline{Case 7}. $c_0=c_1=c_2= n+1$.

By case (iv) and (v) of Definition \ref{main}(1), one of the following
occurs:
\begin{description}
\item{$(\alpha)$} $a_0=a_2=a_3$ and $b_0\unlhd b_1\unlhd b_2$, or
\item{$(\beta)$} $a_0=a_1<_n a_2$ or
\item{$(\gamma)$} $a_0<_n a_1=a_2$ or
\item{$(\delta)$} $a_0<_n a_1 <_n a_2$.
\end{description}
In any case, the conclusion follows from cases (iv) and (v)
of Definition \ref{main}(1).

\underline{Case 8}. $c_0=c_1= n+1$ and $c_2\le n$.

Hence $a_0\le_{n} a_1<_n\langle a_2,b_2,c_2\rangle$ (by cases (iv),
(v) and (iii) of Definition \ref{main}(1)), so the conclusion 
follows by case (v) of Definition \ref{main}(1).

(3) Let $x=
\langle a_0,b_0,c_0\rangle
\neq y=\langle a_1,b_1,c_1\rangle\in C$ be such that $x\le y$. Let $n$
be the first
such that $x,y\in C_n$, hence at least one among $c_0, c_1$ is $n$.
We proceed by induction on $n$.

\underline{$n=0$}. Hence $a_0=a_2=\emptyset$ and $c_0\le c_1=0$,
while $b_0\unlhd b_1$ and $b_0\neq b_1$. For all $s\in {}^{<\omega}2
\setminus\{\langle\rangle\}$ we have
\[
x<\langle x,s,1\rangle <y,
\]
by rules (ii) and (iii) of Definition \ref{main}(1). If $s,t
\in {}^{<\omega}2
\setminus\{\langle\rangle\}$
are incomparable in ${}^{<\omega}2$ we have, by rules (iv) and (v) of
Definition \ref{main}(1), that
\[
\langle x,s,1\rangle><\langle x,t,1\rangle,
\]
hence setting $Z$
to consist of all sequences $\langle x,s,1\rangle$
where $s$ is  a finite sequence of the form $1\dots 10$, 
will satisfy the requirement.

\underline{$n+1$}. Similarly to the discussion
of transitivity, we consider several possibilities. We include
the complete discussion for the sake of completeness.

\smallskip

\underline{Case 1}. $c_0\le n$ and $c_1=n+1$ and $x<_n a_1$.

By the induction hypothesis, there is an infinite
antichain $Z$ of $\langle C,\le\rangle$
such that for all $z\in Z$ we have $x < z< a_1$. By rule (ii)
of Definition \ref{main}(1), we have that $a_1\le y$, so we finish
by transitivity.

\underline{Case 2}. $c_0\le n$ and $c_1=n+1$ and $a_1=x$.

For $s\in {}^{<\omega}2\setminus\{\langle\rangle\}$ consider
$\langle x,s,n+2\rangle$, which is an element of $C_{n+2}$.
For any such $s$ we have
\[
x\le_{n+2}\langle x,s,n+2\rangle\le_{n+2} y,
\]
by rules (ii) and (iii) of Definition \ref{main}(1) respectively.
We finish by transitivity, noticing that for
$s,t\in {}^{<\omega}2\setminus\{\langle\rangle\}$ which are
incomparable in ${}^{<\omega}2$ we have that $\langle x,s,n+2\rangle$
and $\langle x,t,n+2\rangle$ are incomparable in $C$.

\underline{Case 3}. $c_0=n+1$ and $c_1\le n$ and $a_1<_n y$.

For $s\in {}^{<\omega}2 $ such that $b_1\unlhd s$ we have
that $x\le \langle a_1,s,c_1\rangle$, by rule (iv) 
of Definition \ref{main}(1), and we have 
$\langle a_1,s,c_1\rangle\le y$ by rule (iii). Again for $s,t$
incomparable in ${}^{<\omega}2$, we have that 
$\langle a_1,s,c_1\rangle$ and $\langle a_1,t,c_1\rangle$
are incomparable.

\underline{Case 3}. $c_0=n+1=c_1$.

Hence for all $s\in {}^{<\omega}2\setminus\{\langle\rangle\}$
we have that $x\le \langle x,s,n+2\rangle\le y$ by rules
(ii) and (iii) respectively.

\smallskip

Noticing that we have covered all the possibilities, we finish the
proof of the Claim.
$\eop_{\ref{one}}$
\end{Proof of the Claim}

\begin{Claim}\label{two} The set
\[
A\deq\{x\deq\left\langle \langle \emptyset,\langle\rangle,0\rangle,
\langle 0\rangle,1\right\rangle,
y\deq\left\langle \langle \emptyset,\langle\rangle,0\rangle,
\langle 1\rangle,1\right\rangle\}
\]
is a maximal antichain in $\langle C,\le\rangle$.
\end{Claim}

\begin{Proof of the Claim} $A$ is obviously an antichain.
Given $z=\langle a,b,c\rangle\in C$.

If $c=0$, then $a=\emptyset$ and if $b=\langle\rangle$, then
$x<z$ by rule (ii). Otherwise, $x\le z$ by rule (iii).

By induction on $c\ge 1$ we show that $ x\le z$
or $y\le z$ (or both).

Suppose \underline{$c=1$}. If $a=\langle\emptyset,\langle\rangle,0\rangle$,
the conclusion follows by case (iv) of Definition \ref{main}(1).
Otherwise, we use case (v) of Definition \ref{main}(1).

For $\underline{c=n+1}$, we know by the induction hypothesis
that $x\le a$ or $y\le a$. We finish by rule (v).
$\eop_{\ref{two}}$
\end{Proof of the Claim}

\begin{Claim}\label{three} $A$ from Claim \ref{two} does not
split.
\end{Claim}

\begin{Proof of the Claim} As $\left\langle\langle\emptyset,
\langle\rangle,0\rangle,\langle 0,1\rangle,1\right\rangle
\in \UU(x)\setminus [\UU(y)\cup \DD(x)]$
and $\left\langle\langle\emptyset,
\langle\rangle,0\rangle,\langle 1,1\rangle,1\right\rangle
\in \UU(y)\setminus [\UU(x)\cup \DD(y)]$,
and $\langle\emptyset,
\langle\rangle,0\rangle\notin \UU(x)\cup\UU(y)$.
$\eop_{\ref{three}}$
\end{Proof of the Claim}

$\eop_{\ref{example}}$
\end{Proof}

\begin{Remark}\label{gen} The above generalizes to an induction
longer than $\omega$ steps, in the obvious manner.$\eop_{\ref{gen}}$
\end{Remark}

\eject

\begin{sadrzaj}
Parcijalno uredjenje $(P,\le)$ se zove jako gustim, ako u svakom nepraznom
intervalu $(a,b)$ u $P$, postoje dva neuporediva elementa.
Podskup $A$ od $P$ se zove Spernerov sistem, ako su svaka dva
neindenti\v cna elementa u $A$ neuporediva. Ako je $A$ maksimalan,
u odnosu na $\subseteq$, Spernerov sistem, ka\v zemo da je $A$ podijeljen
ako $A=D\bigcup U$ gdje su $D$ i $U$ disjunktni, i svaki element od $P$
je ili $\le$ nego neki element od $D$, ili $\ge$ nego neki element od $U$.
U svom
\v clanku
\cite{AEG}, R. Alshwede, P.L. Erd\"os and N. Graham, su pokazali da svako
kona\v cno jako gusto parcijalno uredjenje, ima osobinu da 
su svi njegovi maksimalni Spernerovi sistemi podijeljeni. U istom
\v clanku autori su dali hipotezu da je isto ta\v cno i za beskona\v cna
jako gusta
prebrojiva uredjenja. U ovom \v clanku, mi smo pokazali da hipoteza
iz \cite{AEG} nije ta\v cna, konstrui\v su\'ci beskona\v cno jako gusto
parcijalno uredjenje koje ima maksimalni Spernerov sistem koji je
nepodijeljen. 
\end{sadrzaj}

\begin{thebibliography}{xxxx}

\bibitem[1]{AEG} R. Ahlswede, P.~L. Erd\"os and N. Graham,
{\em A Splitting Property of Maximal Antichains}, Combinatorica
15 (4) 1995, pp. 475-480.

\bibitem[2]{E} P.~L. Erd\"os, {\em Splitting Property in
Infinite Posets}, Discrete  Mathematics, 163 (1-3) 1997, pp. 251-256.

\end{thebibliography}
\end{document}